\documentclass[12pt]{article}

\usepackage{amsfonts}

\newtheorem{thm}{Theorem}[section]
\newtheorem{lemma}[thm]{Lemma}
\newtheorem{cor}[thm]{Corollary}

\newtheorem{prop}[thm]{Proposition}
\newtheorem{conjecture}{Conjecture}

\newcommand{\be}{\begin{equation}}
\newcommand{\ee}{\end{equation}}
\newcommand{\openbox}{\leavevmode
  \hbox to8pt{\hfil\vrule\vbox to6pt{\hrule width6pt\vfil\hrule}\vrule}}

\newcommand{\qed}{\hbox to5pt{ } \hfill \openbox\bigskip\medskip}

\newcommand{\ve}[1]{\mathbf{#1}}
\newcommand{\cT}{\mbox{$\cal T$}}

\newcommand{\cF}{\mbox{$\cal F$}}

\newcommand{\cH}{\mbox{$\cal H$}}

\newcommand{\cR}{\mbox{$\cal R$}}

\newcommand{\cQ}{\mbox{$\cal Q$}}
\newcommand{\cO}{\mbox{$\cal O$}}

\newcommand{\N}{\mathbb N}

\newcommand{\R}{\mathbb R}

\newcommand{\F}{\mathbb F}

\title{Upper bounds for the size of set systems with a symmetric  set of Hamming distances}
\author{G\'abor Heged\"{u}s
\\{\normalsize  \'Obuda University}
\\{\normalsize B\'ecsi \'ut 96/B, Budapest, Hungary, H-1032}
\\{\normalsize hegedus.gabor@uni-obuda.hu}
}

\begin{document}

\maketitle
\begin{abstract}
Let $\mbox{$\cal F$}\subseteq 2^{[n]}$ be a fixed family of subsets. Let $D(\mbox{$\cal F$})$ stand for  the following set of Hamming distances: 
$$
D(\mbox{$\cal F$}):=\{d_H(F,G):~ F, G\in \mbox{$\cal F$},\ F\neq G\}.
$$  
$\mbox{$\cal F$}$ is said to be a  Hamming symmetric family, if $d\in D(\mbox{$\cal F$})$ implies $n-d\in D(\mbox{$\cal F$})$ for each $d\in D(\mbox{$\cal F$})$.

We give sharp upper bounds for the size of Hamming symmetric families. Our proof is based on  the linear algebra bound method.
\end{abstract}

\medskip

\section{Introduction}

Throughout the paper 
$n$ denotes a positive integer and $[n]$ stands for the set $\{1,2,
\ldots, n\}$. The family of all subsets of $[n]$ is denoted by $2^{[n]}$. 
For an integer $0\leq d\leq n$ we denote by
${[n] \choose d}$ the family of all  $d$ element subsets of $[n]$,
 and ${[n] \choose \leq d}={[n] \choose 0}\cup\ldots\cup{[n] \choose d}$
the subsets of size at most $d$.

We say that $\mbox{$\cal F$}$ is a {\em complete intersecting family}, if $\cF$ is an intersecting family and $|\cF|=2^{n-1}$. For example, the family $\cH:=\{G\subseteq  [n]:~ 1\in G\}$ is  a  complete intersecting family.

Let $\cF\subseteq 2^{[n]}$ be a fixed family of subsets and $F, G\in \cF$ be two distinct elements of $\cF$.
Let $d_H(F,G)$ stand for the Hamming distance of the sets
 $F$ and  $G$, i.e.,  
$d_H(F,G):=|F \Delta G|$, 
where $F \Delta G$ is the usual symmetric difference.

Denote by $D(\cF)$ the following set of Hamming distances: 
$$
D(\cF):=\{d_H(F,G):~ F, G\in \cF,\ F\neq G\}.
$$  
       
Let $q>1$ be an integer. Let $\cH\subseteq \{0,1,\ldots ,q-1\}^n$ and 
let $\ve h_1,\ve h_2\in \cH$ be two elements of the vector system $\cH$.  
Let $d_H(\ve h_1,\ve h_2)$ stand for the Hamming distance of the
vectors $\ve h_1,\ve h_2\in \cH$:
$$
d_H(\ve h_1,\ve h_2):=|\{i\in [n]:~  (\ve h_1)_i\neq (\ve h_2)_i\}|.
$$

Denote by $D(\cH)$ the following set of Hamming distances: 
$$
D(\cH):=\{ d_H(\ve h_1,\ve h_2):~ \ve h_1,\ve h_2\in \cH, \ve h_1\neq \ve h_2 \}.
$$  

Delsarte proved the following well-known upper bound for the size of the vectors systems with $s$ distinct Hamming distances (see in \cite{D},  \cite{D2}).

\begin{thm} \label{Delsarte}
Let $0<s\leq n$, $q>1$ be positive integers. 
 Let $L=\{\ell_1 ,\ldots ,\ell_s\}\subseteq [n]$ be a set of $s$ positive integers. Let $\cH\subseteq \{0,1,\ldots ,q-1\}^n$ and suppose that $d_H(\ve h_1,\ve h_2)\in L$ for each distinct $\ve h_1,\ve h_2\in \cH$ vectors. Then
$$
|\cH|\leq \sum_{i=0}^s {n\choose i} (q-1)^i.
$$
\end{thm}

In the $q=2$ special case we get the following statement.

\begin{cor} \label{Delsarte2}
Let $0<s\leq n$ be positive integers. 
Let $\cF\subseteq 2^{[n]}$ be a set system such that $|D(\cF)|=s$. Then
$$
|\cF|\leq \sum_{j=0}^s {n\choose j}.
$$
\end{cor}

Let $\cF$ be  a family of subsets of $[n]$.
$\cF$ is said to be a  {\em Hamming symmetric family}, if $d\in D(\cF)$ implies $n-d\in D(\cF)$ for each $d\in D(\cF)$. Specially  if $\cF$ is a Hamming symmetric family, then $n\not\in D(\cF)$.

Our main result follows.

\begin{thm} \label{main}
Let $\cF$ be  a Hamming symmetric family of subsets of $[n]$. Let $s:=|D(\cF)|$. 

If $n/2\notin D(\cF)$, then 
$$
|\cF|\leq \sum_{j=0}^{\lfloor \frac s2 \rfloor} {n\choose 2j}.
$$
If $n/2\in D(\cF)$, then
$$
|\cF|\leq \sum_{j=0}^{\lfloor \frac{s-1}{2} \rfloor} {n\choose 2j+1}.
$$
\end{thm}

It is easy to verify that Theorem \ref{main} is sharp. Let $n=2t$ be an even integer and consider the family
$$
\cF:=\{G\subseteq  [n]:~ 1\in G\}.
$$ 
Then $\cF$ a complete intersecting family with  $|\cF|=2^{n-1}$. 
Clearly $D(\cF)=[n-1]$ and $s=n-1$.  Hence $n/2\in D(\cF)$ and 
$$
|\cF|=\sum_{j=0}^{\lfloor \frac{n-2}{2} \rfloor} {n\choose 2j+1}=2^{n-1}.
$$

\section{Preliminaries}

The proof of  our main result is based on  the linear algebra bound method and the Determinant Criterion (see \cite{BF} Proposition 2.7). First we recall here shortly for the reader's convenience this principle.

\begin{prop} \label{det} (Determinant Criterion)
Let $\F$ denote an arbitrary field. Let $f_i:\Omega \to \F$ be functions and $\ve v_j\in \Omega$ elements for each $1\leq i,j\leq m$  such that the $m \times m$ matrix $B=(f_i(\ve v_j))_{i,j=1}^m$
is non-singular. Then $f_1,\ldots ,f_m$ are linearly independent functions of the space $\F^{\Omega}$. 
\end{prop}

Let $0\leq s\leq n$, $q>1$ be integers.

We define the following sets of monomials:
$$
\cT(n,s):=\{\alpha=(\alpha_1, \ldots , \alpha_n)\in {\N}^n:~  \sum_{i=1}^n \alpha_i \mbox{ is even },\ \sum_{i=1}^n\alpha_i\leq s\},
$$
$$
\cQ(n,s):=\{\alpha=(\alpha_1, \ldots , \alpha_n)\in \cT(n,s):~ \alpha_i\leq 1 \mbox{ for each }1\leq i\leq n\},
$$
$$
\cO(n,s):=\{\alpha=(\alpha_1, \ldots , \alpha_n)\in {\N}^n:~  \sum_{i=1}^n \alpha_i \mbox{ is odd },\ \sum_{i=1}^n\alpha_i\leq s\}
$$
and
$$
\cR(n,s):=\{\alpha=(\alpha_1, \ldots , \alpha_n)\in \cO(n,s):~ \alpha_i\leq 1 \mbox{ for each }1\leq i\leq n\}.
$$

We use the following combinatorial lemmas in the proofs of our main results.
\begin{lemma} \label{meret}
Let $0\leq s\leq n$ be integers. Then
$$
|Q(n,s)|=\sum_{j=0}^{\lfloor \frac s2 \rfloor} {n\choose 2j}.
$$
\end{lemma}
\begin{lemma} \label{meret2}
Let $0\leq s\leq n$ be integers. Then
$$
|R(n,s)|=\sum_{j=0}^{\lfloor \frac{s-1}{2} \rfloor} {n\choose 2j+1}.
$$
\end{lemma}

\section{Proof}

Consider the set
$$
D'(\cF):=\{n-2d:~ d\in D(\cF)\}.
$$

Let $\ve v_F\in \{-1,1\}^n$ denote the characteristic vector of the set $F$, i.e. 
$\ve v_F(k)=1$, if $k\in F$, $\ve v_F(k)=-1$ otherwise.

If $F,G\in 2^{[n]}$ are arbitrary subsets , then it is easy to check that
$$
n-2d_H(F,G)=\langle \ve v_F,\ve v_G\rangle. 
$$
It follows from this equality that  $D'(\cF)$ is precisely the following set of scalar products:
$$
D'(\cF):=\{\langle \ve v_F,\ve v_G\rangle:~ F, G\in \cF,\ F\neq G\}.
$$
It is clear that if $n/2\notin D(\cF)$, then $0\notin D'(\cF)$.

{\bf Proof of Theorem \ref{main}:}

Consider the real polynomials
$$
g_j(x_1,\ldots ,x_n)=x_j^2 -1\in \R[x_1, \ldots ,x_n]
$$
for each  $1\leq j\leq n$.

Let $\cF=\{F_1, \ldots ,F_m\}$ be  a Hamming symmetric family of subsets of $[n]$, where  $s=|D(\cF)|$. Then  it is easy to see that $D'(\cF)$ is symmetric set with respect to $0$, i.e., if $d\in D'(\cF)$, then 
 $-d\in D'(\cF)$.

Let $i\in [m]$ be a fixed index. Let $\ve v_i\in \{-1,1\}^n$ denote the characteristic vector of the set $F_i$.

1. First suppose that $n/2\notin D(\cF)$. Then  $0\notin  D'(\cF)$.



Consider the polynomials
\begin{equation}  \label{Ppol} 
P_i(x_1, \ldots ,x_n):= \prod_{d\in D'(\cF)} \Big( \langle \ve x, \ve v_i\rangle-d \Big)\in \R[\ve x]
\end{equation}
for each $1\leq i\leq m$, where $\langle \ve x, \ve v_i\rangle$ denotes the usual scalar product of the vectors $\ve x$ and $ \ve v_i$. Clearly $\mbox{deg}(P_i)\leq s=|D'(\cF)|$ for each $1\leq i\leq m$.

$D'(\cF)$ is symmetric with respect to $0$, consequently if we expand $P_i$  as a linear combination of monomials, then we get 
\begin{equation}  \label{exp} 
P_i(x_1, \ldots ,x_n)=\sum_{\alpha\in \cT(n,s)} d_{\alpha}x^{\alpha},
\end{equation}
where $d_{\alpha}\in \R$ are real coefficients  for each ${\alpha}\in \cT(n,s)$ (here $x^{\alpha}$ denotes the monomial  $x_{1}^{\alpha_{1}}\cdot \ldots \cdot x_{n}^{\alpha_{n}}$). 

But $\ve v_i\in \{-1,1\}^n$ for each $i$,  consequently the equation 
\begin{equation}  \label{kocka} 
x_j^2=1
\end{equation}
is true for each $\ve v_i$ vector and  for each $1\leq j\leq n$. 

Let $Q_i$ denote the polynomial obtained by writing  $P_i$ as a linear combination of monomials and replacing,
repeatedly, each occurrence of $x_j^2$, where $1\leq j\leq n$,  by $1$.

Since $g_t(\ve v_i)=0$ for each $i$ and  for each $1\leq t\leq n$, 
hence $Q_i(\ve v_j)=P_i(\ve v_j)$ for each $1\leq i\ne j\leq r$. 

We prove that the set of polynomials  $\{Q_i:~ 1\leq i\leq m\}$ is linearly independent. This fact follows from the Determinant Criterion, when we define $\F:=\R$, $\Omega= \{-1,1\}^n$ and $f_i:=Q_i$ for each $i$.  It is enough to prove that  $Q_i(\ve v_i)=P_i(\ve v_i)\ne 0$ for each $1\leq i\leq m$ and $Q_i(\ve v_j)=P_i(\ve v_j)=0$ for each $1\leq i\ne j\leq m$, since then we can apply the Determinant Criterion.

But $P_i(\ve v_i)=\prod_{d\in D'(\cF)} (n-d)$ for each $1\leq i\leq m$.  Consequently  $P_i(\ve v_i)\ne 0$ for each $1\leq i\leq m$, since $n\notin D'(\cF)$.

It follows from the definition of $D'(\cF)$ that $P_i(\ve v_j)=\prod_{d\in D'(\cF)} \Big((\langle \ve v_i, \ve v_j\rangle)-d \Big)=0$ for each $1\leq i\ne j\leq m$.

Then it is easy to verify that  we can write $Q_i$ as a linear combination of monomials in the form
$$
Q_i=\sum _{\alpha\in \cQ(n,s)} t_{\alpha}x^{\alpha},
$$ 
where $t_{\alpha}\in \R$ are the real coefficients for each ${\alpha}\in \cQ(n,s)$. This follows immediately from the expansion (\ref{exp}) and from the relations (\ref{kocka}).  

Since the polynomials $\{Q_i:~ 1\leq i\leq m\}$ are linearly independent and 
if we expand $Q_i$ as a linear combination of monomials, then
all monomials appearing in this linear combination  contained in the set of monomials 
$$
\{x^{\alpha}:~ \alpha\in \cQ(n,s)\}
$$
for each $i$,  we infer from Lemma  \ref{meret} that 
$$
m=|\cF|\leq  |\cQ(n,s)|=\sum_{j=0}^{\lfloor \frac s2 \rfloor} {n\choose 2j}.
$$
\qed.

2. Suppose that $n/2\in D(\cF)$. Then  $0\in  D'(\cF)$.

We can define the $P_i$ polynomials precisely the same way as in the first part of the proof (see the definition appearing in (\ref{Ppol}). 

But if we expand $P_i$  as a linear combination of monomials, then we get 
\begin{equation}  \label{exp2} 
P_i(x_1, \ldots ,x_n)=\sum_{\alpha\in \cO(n,s)} c_{\alpha}x^{\alpha},
\end{equation}
where $c_{\alpha}\in \R$ are real coefficients  for each ${\alpha}\in \cO(n,s)$. 

If we reduce the polynomial  $P_i$ with the equations (\ref{kocka}), then we get a new polynomial $R_i$.
Clearly  $R_i(\ve v_j)=P_i(\ve v_j)$ for each $1\leq i\ne j\leq r$.

It can be proved easily from the Determinant Criterion that the set of polynomials  $\{R_i:~ 1\leq i\leq m\}$ is linearly independent.

Since if we expand $R_i$ as a linear combination of monomials, then
all monomials appearing in this linear combination  contained in the set of monomials 
$$
\{x^{\alpha}:~ \alpha\in \cR(n,s)\}
$$
for each $i$, hence we get that
$$
m=|\cF|\leq  |\cR(n,s)|=\sum_{j=0}^{\lfloor \frac{s-1}{2} \rfloor} {n\choose 2j+1}
$$
by Lemma \ref{meret2}.

\section{Concluding remarks}

We propose the following conjecture as a generalization of Theorem \ref{main}.

\begin{conjecture} \label{Hconj}
Let $0<s\leq n$, $q>1$ be positive integers. 
 Let $L=\{\ell_1 ,\ldots ,\ell_s\}\subseteq [n]$ be a set of $s$ positive integers. Let $\cH\subseteq \{0,1,\ldots ,q-1\}^n$ 
Let $\cH\subseteq \{0,1,\ldots ,q-1\}^n$ be  a Hamming symmetric vector system such that $d_H(\ve h_1,\ve h_2)\in L$ for each distinct $\ve h_1,\ve h_2\in \cH$ vectors.

If $n/2\notin D(\cF)$, then 
$$
|\cH|\leq \sum_{0\leq  i\leq s,i\equiv 0 \pmod 2} {n\choose i} (q-1)^i.
$$
If $n/2\in D(\cF)$, then
$$
|\cF|\leq \sum_{0\leq  i\leq s,i\equiv 1 \pmod 2} {n\choose i} (q-1)^i.
$$
\end{conjecture}


\end{document}